\documentclass[12pt]{article}
\usepackage{bbm}
\usepackage{pifont}
\usepackage{bbding}
\usepackage{stmaryrd}
\usepackage{amssymb}
\usepackage{amsthm}
\usepackage{fancyhdr}
\usepackage{amsfonts}
\usepackage{mathrsfs}
\usepackage{graphics}
\usepackage{graphicx}
\usepackage{amsmath}

\newtheorem{theorem}{Theorem}

\newtheorem{Open Problem}[theorem]{Open Problem}

\usepackage{amsthm}
\newtheorem*{Open}{Open Problem}
\newtheorem*{Conj}{Conjecture}
\begin{document}
\baselineskip=0.7truecm
\begin{center}\huge\Large The maximum forcing number of polyomino\footnote{Supported by NSFC
(Grant No.11301217,11171134), the Fundamental
Research Funds for the Central Universities (Grant No.2010121007),
the Natural Science Foundation of Fujian Province (Grant
No.2013J01014) and New Century Excellent
Talents in Fujian Province University(grant JA14168).\\Corresponding author: xuliqiong@jmu.edu.cn}
 \normalsize \\ [6truemm] Yuqing Lin$^{a}$, Liqiong Xu$^{b}$ and Fuji Zhang$^{c}$\\
 {\small$^a$ School of Electrical Engineering and Computer Science, University of Newcastle,}\\{\small NSW2308, Australia
}\\
 {\small $^b$ School of Mathematical Sciences, Jimei University,}\\
{\small Jimei, Fujian 361005,
 P. R. China}\\
 {\small$^c$ School of Mathematical Sciences, Xiamen University,}\\{\small Xiamen, Fujian 361005,
 P. R. China}
\end{center}

\normalsize\thispagestyle{empty}

\noindent{\sl\bf Abtract}\ \ The forcing number of a perfect
matching $M$ of a graph $G$ is the cardinality of the smallest
subset of $M$ that is contained in no other perfect matchings of $G$.
For a planar embedding of a 2-connected bipartite planar graph $G$ which has a perfect matching,
the concept
of Clar number of hexagonal system had been extended by Abeledo and Atkinson as follows:
a spanning subgraph $C$ of  is called a Clar cover of $G$ if each of its components is either an even face or an edge, the maximum number of even faces in Clar covers of $G$ is
called Clar number of $G$, and the Clar cover with the maximum number of even faces is called the maximum Clar cover.
It was proved that if $G$ is a hexagonal system with a perfect matching $M$ and $K'$ is a
set of hexagons in a maximum Clar cover of $G$, then $G-K'$ has a unique
1-factor. Using this result, Xu {\it et. at.} proved that the maximum forcing number of the elementary hexagonal system are equal to their Clar numbers, and then the maximum forcing
number of the elementary hexagonal system can be computed in polynomial time. In this paper, we show that an elementary
polyomino has a unique perfect matching when removing the set of tetragons from its maximum Clar cover. Thus the maximum forcing number of elementary polyomino equals to its Clar number and can be computed in polynomial time. Also, we have extended our result to the non-elementary polyomino and hexagonal system.

\section{Introduction}
\hspace*{8mm} Let $G$ be a graph that has a perfect matching. A
forcing set for a perfect matching $M$ of $G$ is a subset $S$ of
$M$, such that $S$ is contained in no other perfect matchings of $G$.
The cardinality of a smallest forcing set of $M$ is
called the forcing number of $M$, and is denoted by $f(G;M)$. The
minimum and maximum of $f(G;M)$ over all perfect
matchings $M$ of $G$ is denoted by $f(G)$ and $F(G)$, respectively.
Given a matching $M$ in a graph $G$, an $M$-alternating path (cycle) is
a path (cycle) in $G$ whose edges are alternately in $M$ and outside
of $M$. Let $e$ be an edge of $G$. If $e$ is contained in all perfect
matchings of $G$, or is not contained in any perfect matchings of $G$, then $e$ is called fixed
bond. If $e$ is contained in a perfect
matchings of $G$, then $e$ is called double
bond.

A polyomino (respectively, hexagonal system) is a finite connected plane graph with no cut
vertex and every interior region is surrounded by a regular square
(respectively, hexagon). A connected bipartite graph is called elementary (or normal) if its
every edge is contained in some perfect matching. Let $G$ be a plane
bipartite graph, a face of $G$ is called resonant if its boundary is
an alternating cycle with respect to a perfect matching of $G$.

The Clar number is originally defined for hexagonal systems by Clar[1]. The Clar cover polynomial had been introduced in [2,3]. The further progress of Clar's theory can be found in two surveys [4,5].
 Later, Abeledo and Atkinson[6] generalized the concept of Clar number
for bipartite and 2-connected plane graphs. For a planar embedding of a 2-connected bipartite planar $G$,
a Clar cover of $G$ is a spanning subgraph $C$ with its components is either an even face or
an edge, the maximum number of even faces in Clar covers of $G$ is
called Clar number of $G$, and denoted by $C(G)$. We call a Clar cover with the maximum number of even face the maximum Clar cover.

The idea of forcing number was inspired by practical chemistry
problems.  This concept was first proposed by Harary {\it et.
al.} in [7]. The same idea appeared in earlier papers by
Randi{\'c} and Klein [8,9] in terms of ¡°innate
degree of freedom¡± of a Kekul{\'e} structure. The forcing
numbers of  square grids, stop
signs, torus, hypercube, and some special fullerene
graphs have been considered in [10-16].
 Adams et al.\cite{Adams} proved it is
NP-complete to find the smallest forcing set of a bipartite graph
with maximum degree 3. Later, Afshani etal.[10] proved it is
NP-complete to find the smallest forcing number of a bipartite graph
with maximum degree 4, and proposed the following question:

\begin{Open} Given a
graph $G$, what is the computational complexity of finding the
maximum forcing number of $G$.
\end{Open}

Recently, Xu,Bian and Zhang\cite{Xu} showed that when $G$ is an elementary hexagonal system,
the maximum forcing number of $G$ can be computed in polynomial time. In this paper,  we first show that an elementary
polyomino has a unique perfect matching when removing all tetragons of a maximum Clar cover of the polyomino. Base on this result, we know that the maximum
forcing number of an elementary polyomino equals to its Clar
number, thus the  maximum forcing number of the elementary polyomino  can be computed in polynomial time.
Finally, we also show that for both elementary and non-elementary hexagonal system and polyomino, their maximum forcing number can be computed in polynomial time, which generalize the result in \cite{Xu} to non-elementary hexagonal system.

\section{A maximum Clar cover of polyomino}

\ \ \ \ \ \ \ First we give some theorems in the following.

\begin{theorem}[{\cite{Xu}}] Let $G$ be a plane elementary bipartite graph. Then $F(G) \geq C(G)$.\end{theorem}

\begin{theorem}[{\cite{Mao}}]\label{Thm-2.7} Let $H$ be a hexagonal system with a perfect matching $M$ and $C'$ be a
set of hexagons in a Clar cover of $H$ with maximum cardinality of hexagons. Then $H-C'$ has a unique
1-factor.
\end{theorem}

In \cite{Hansen}, Hansen and Zheng formulate the computation of Clar number as an
 integer program problem. Later, Abeledo and Atkinson\cite{Abeledo}  proved that:

\begin{theorem}[{\cite{Abeledo}}]\label{Thm-2.7} Let $G$ be a 2-connected bipartite planar graph. Then the Clar number of $G$ can be
computed in polynomial time using linear programming methods.
\end{theorem}

Zhang et. al. showed that

\begin{theorem}[{\cite{Zhang5}}]\label{Thm-2.5}Let $G$ be a plane
bipartite graph with more than two vertices. Then each face of $G$
is resonant if and only if $G$ is elementary.
\end{theorem}

Using Theorems 1,2,3 and 4,  Xu,Bian and Zhang \cite{Xu} proved that the maximum forcing number of an elementary hexagonal system equals to its Clar number, and then the maximum forcing number of the elementary hexagonal system can be computed in polynomial time, and conjectured that

\begin{Conj} {\cite{Xu}} Let $G$ be an elementary polyomino, then the maximum forcing number of $G$ can be computed in
polynomial time.
\end{Conj}

We now prove that the conjecture is true.

\begin{theorem}\label{Thm-2.8} Let $G$ be an elementary polyomino with perfect matchings. Let $K$ be a maximum Clar cover of $G$ with $C(G)$ tetragons, and
let $K'$ be the set of tetragons in $K$. Then $G-K'$ has a unique 1-factor.
\end{theorem}

{\bf Proof}: Suppose that $G - K'$ has more than one perfect matching. Joining any two of the perfect matchings of $G - K'$ will give us a set of alternating cycles in the graph $G$. Choose the two of the perfect matchings of $G - K'$ such that one of the alternating cycles $C$ is smallest possible. Let $G^*$ denote the subgraph of $G$ such that the outer boundary of $G^*$ is the alternating cycle $C$. Clearly, there is no alternating cycles for any perfect matchings in $G^*$ since $C$ is the smallest.

Now let's look at the graph $G^*$ and we label the graph $G^*$ in the following way:

The top left vertex is labelled $v_{1,1}$, and the vertices in the same row will be labelled as $v_{1,i}$. The second row vertices is labelled as $v_{2,i}$ and $i$ could be 0 or a negative value if the vertex is on the left of the vertex $v_{2,1}$. The face with vertex $v_{i,j}$ on its left top is labelled as $f_{i,j}$, see Fig.1 for an example. Let $K'|_{G^*}$ denote the set of tetragons of $K'$ in $G^*$. By the choose of perfect matchings, $C$ is a nice cycle of $G$(i.e. $G-C$ has a perfect matching) and $G^* -C-K'|_{G^*}$ has a unique perfect matching which is denoted by $M(G^*)$.

{\bf Observation 1}
In $M(G^*)$, if there are two parallel edges, i,e, two opposite edges of a face, either $v_{i,j}v_{i,j+1}$ and $v_{i+1,j}v_{i+1,j+1}$ or $v_{i,j}v_{i+1,j}$ and $v_{i,j+1}v_{i+1,j+1}$, then $K$ is not a maximum Clar cover.

Clearly, then the tetragon surrounded by these two parallel edges can be added to $K'$ and then $K$ is not a maximum Clar cover.

Along the same line of reasoning, we have

{\bf Observation 2}
If there are two parallel edges on the boundary of $G^*$, i,e, edges such as $v_{i,j}v_{i,j+1}$ and $v_{i+1,j}v_{i+1,j+1}$ or $v_{i,j}v_{i+1,j}$ and $v_{i,j+1}v_{i+1,j+1}$ on $C$,
then $K$ is not a maximum Clar cover.

Suppose that the face $f_{i,j}$ is surrounded by these two parallel edges, clearly that if $f_{i,j}$ is a pending face, then the corresponding tetragon can be included in the $K'$, a contradiction. In other cases, removing the tetragon that corresponding to face $f_{i,j}$, the $M(G^*)$ has not been disturbed. Furthermore, the cycle $C$ has been disconnected into two odd length paths with a unique perfect matching. Thus, we could include the tetragon corresponding to face $f_{i,j}$ to increase $K'$, a contradiction.

{\bf Observation 3}
In $M(G^*)$, If there is an edge parallel to the boundary, then $K$ is not a maximum Clar cover.


Now we shall prove that either there exist a pair of parallel edges in $K(G^*)$ or we could replace some tetragons of $K'$ by a larger set, which lead to a contradiction.

Look at the consecutive faces $f_{1,1}, f_{1,2},...f_{1,i}$ in the first row of $G^*$. Because of Observation 2, we know that $i > 1$. Now we look at the faces $f_{2,1}, f_{2,2},...f_{2,i}$. Let's first see that these faces are all belong to $G^*$. Suppose that one of the faces $f_{2,j}$, where $1 \leq j \leq i$, is not part of the $G^*$, then $G^*$ has two parallel edges on the boundary of $G^*$. Based on Observation 2, we know it is not possible (See Fig.2 for an example).

\vskip 2mm
\begin{figure}[htbp]
\centering
\includegraphics[height=4cm]{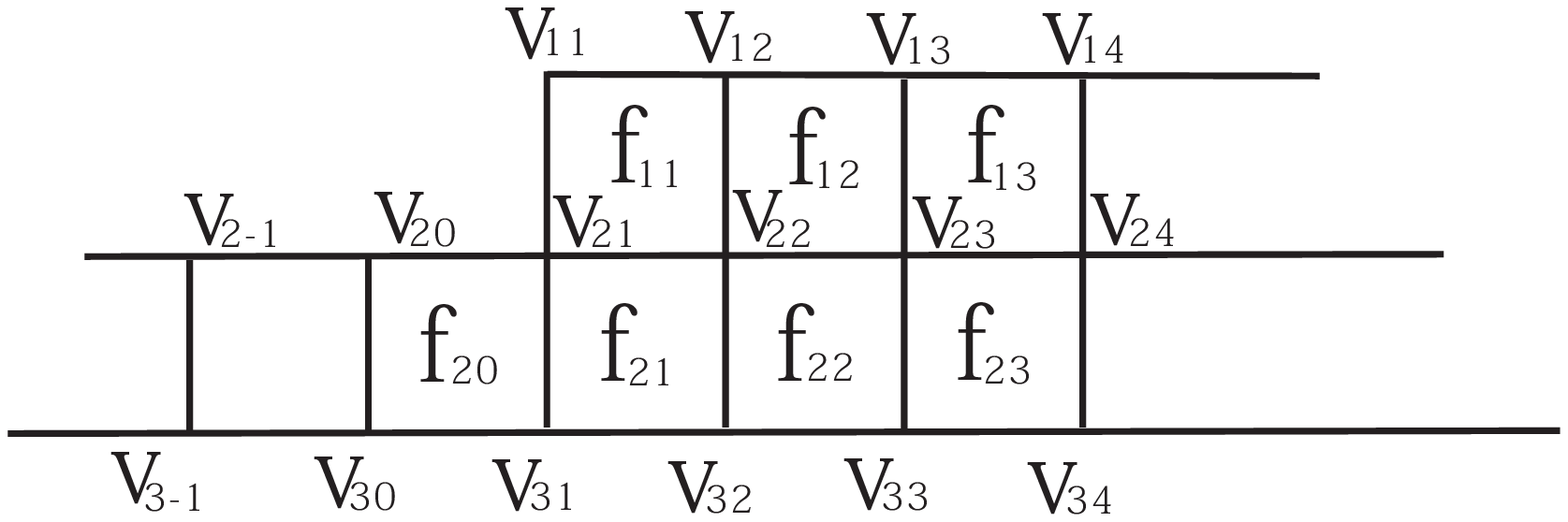}
\par {\small Fig.1.}
\end{figure}

Now, we know that  $f_{2,1}, f_{2,2},...f_{2,i}$ are faces of $G^*$. First, it is clear that none of the edges $v_{2,1}v_{2,2}$, $v_{2,2}v_{2,3}$ .. $v_{2,i-1}v_{2,i}$ are in the $M(G^*)$ because of Observation 3. This implies that the edges $v_{2,2}v_{3,2}$, $v_{2,3}v_{3,3}$. $v_{2,i}v_{3,i}$ are either in the perfect matching $M(G^*)$ or in $K'_{G^*}$. And furthermore, due to Observation 2, there are no parallel edges in the $M(G^*)$.

Suppose that $i$ is odd, then we know that faces $f_{2,2}, f_{2,4},...f_{2,i-1}$ should be in $K'_{G^*}$, otherwise, remove those faces in  $K'$ of form $f_{2,t}$ where $1 \leq t \leq i-1$ and then take $f_{2,2}, f_{2,4},...f_{2,i-1}$ with the remaining tetragons of $K'$ to get a Clar cover of $G$ with larger number of tetragons, a contradiction.

Now we have $f_{2,2}, f_{2,4},...f_{2,i-1} \in K'_{G^*}$, we shall replace $f_{2,2}, f_{2,4},...f_{2,i-1}$ with the faces $f_{1,1}, f_{1,3},...f_{1,i}$ which give us a Clar cover with one more tetragons than that in $K$ (See Fig.3 for an example), a contradiction.

\vskip 2mm
\begin{figure}[htbp]
\centering
\includegraphics[height=3cm]{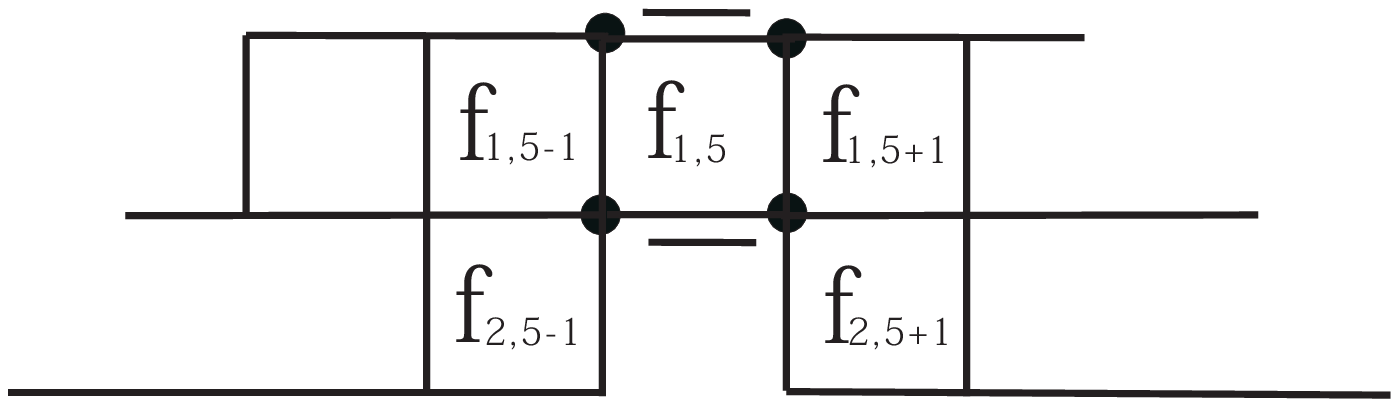}
\par {\small Fig.2.}
\end{figure}

Now we look at the $i$ even case, where the previous suggested approach doesn't work since it will not give us a new Clar with larger number of tetragons. However, we could assume that the edge $v_{2,2}v_{3,2}$ is in $M(G^*)$ and
faces $f_{2,3}, f_{2,5},...f_{2,i-1}$ are in $C'$, otherwise, we could make the rearrangement for that configuration to happen. If the edge $v_{2,1}v_{3,1}$ is on the boundary of $G^*$, then by the Observation 3, we know that $K$ is not a maximum Clar cover of $G$. We then could assume that there are other faces on the left of $f_{2,1}$

Next, we know that the edge $v_{3,1}v_{4,1}$ is either in $M(G^*)$ or $K'_{G^*}$ or belong to the boundary of $G^*$. First we see it is not possible for $v_{3,1}v_{4,1}$ to be on the boundary, since it implies that the edges $v_{2,0}v_{2,1}$ $v_{3,0}v_{3,1}$ are on the cycle $C$, based on the Observation 2, it is a contradiction. Thus we know that $v_{3,1}v_{4,1}$ must belong to $K(G^*)$ or $C'|_{G^*}$.

Assume that the left most vertex on the second row of $G^*$ is $v_{2,-j}$, i.e. it doesn't exist a vertex $v_{2,-t}$ where $t>j$. We could assume that $f_{2,0}, f_{2,-1},..., f_{2,-j}$ are all belong to $G^*$, otherwise, we re-label the graph, take the left most $v_{2,-j}$ as $v_{1,1}$.

Now we also know that all faces $f_{3,0}, f_{3,-1},..., f_{3,-j}$ are belong to $G^*$, otherwise, based on the Observation 2, we could show that one of the faces $f_{2,-t}$, where $1<t<j$, could be included in  $K'_{G^*}$. See Fig.4  for detail.

\vskip 2mm
\begin{figure}[htbp]
\centering
\includegraphics[height=3cm]{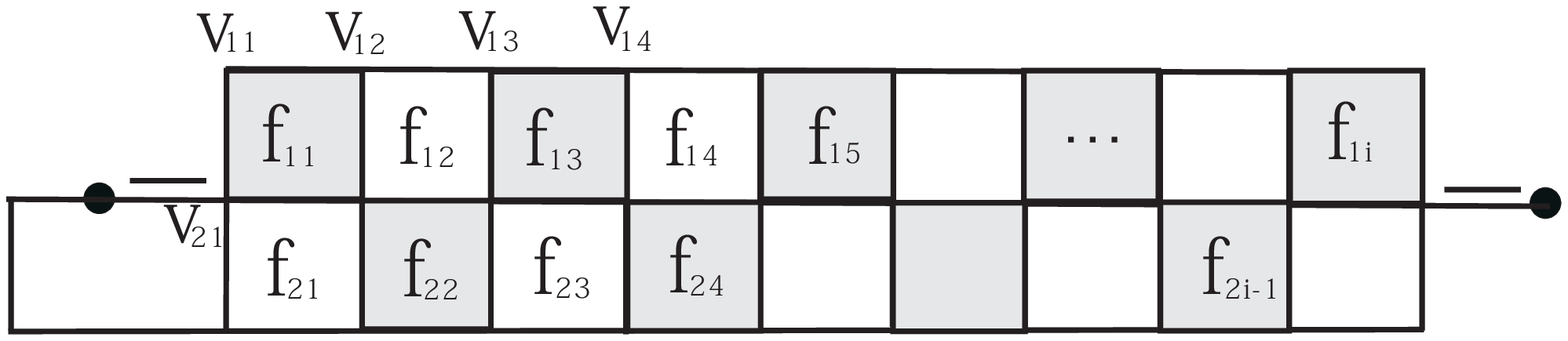}
\par {\small Fig.3.}
\end{figure}

As in previous cases, we know that the edges  $v_{3,0}v_{4,0}$, $v_{3,-1}v_{4,-1}$ .. $v_{3,-j+1}v_{4,-j+1}$ are either in the perfect matching $M(G^*)$ or in $K'_{G^*}$. If $j$ is odd, then faces $f_{3,0}, f_{3,-2},...f_{3,-j+1}$ must be in $C'|_{G^*}$.

If $j$ is odd, we could then replace these faces $f_{3,0}, f_{3,-2},...f_{3,-j+1}$ by $f_{2,1}, f_{2,-1},...f_{2,-j}$, and clearly we have a larger Clar cover of $G$, see Fig.5 for an example. Clearly the left over graph has a perfect matching, i.e. there are odd length of path been removed from the boundary and no internal matchings been disturbed.

\vskip 2mm
\begin{figure}[htbp]
\centering
\includegraphics[height=4cm]{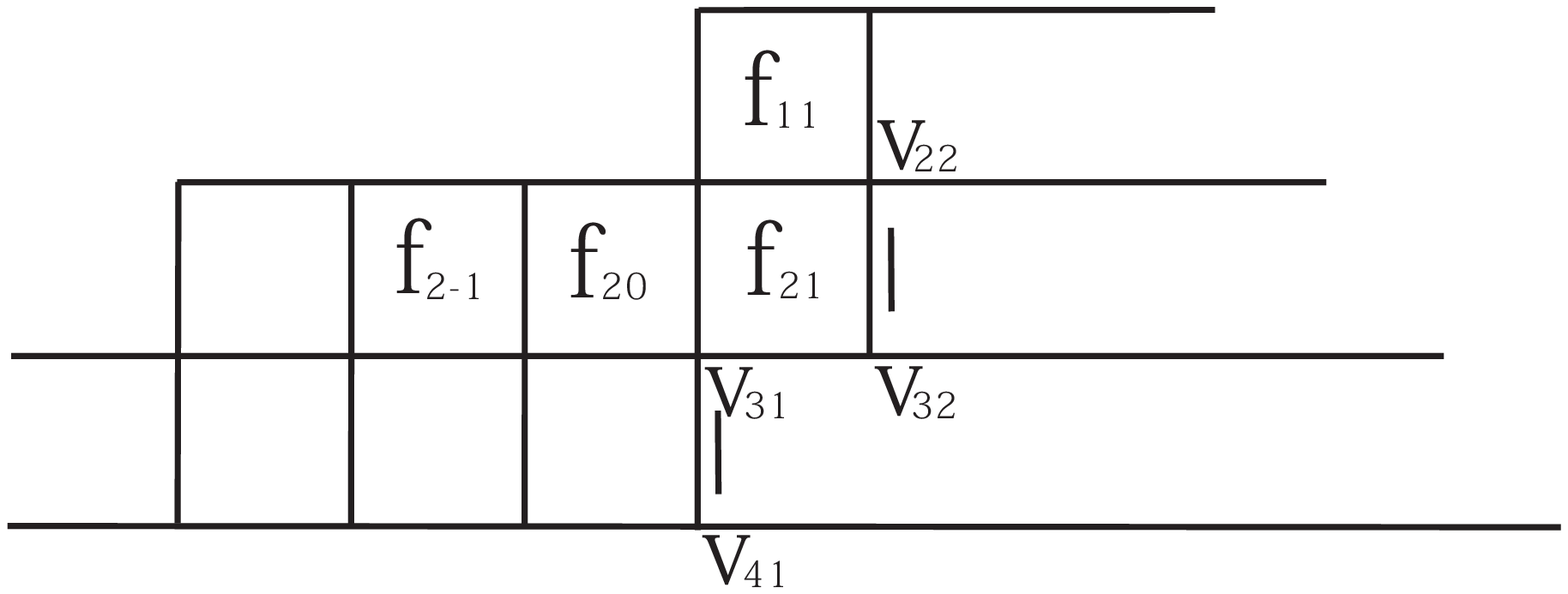}
\par {\small Fig.4.}
\end{figure}

If $j$ is even, then we could assume that the edge $v_{3,-j+1}v_{4,-j+1}$ is in the perfect matching $M(G^*)$, if $v_{3,-j}v_{4,-j}$ is on the boundary of $G^*$, then based on Observation 2, we know $K'(G^*)$ is not maximum. See Fig 6 for an example, now the left over case is that there are more faces on the left of $f_{3,-j}$. Suppose that the left most faces in the third row is $f_{3,-t}$.
\vskip 2mm
\begin{figure}[htbp]
\centering
\includegraphics[height=5cm]{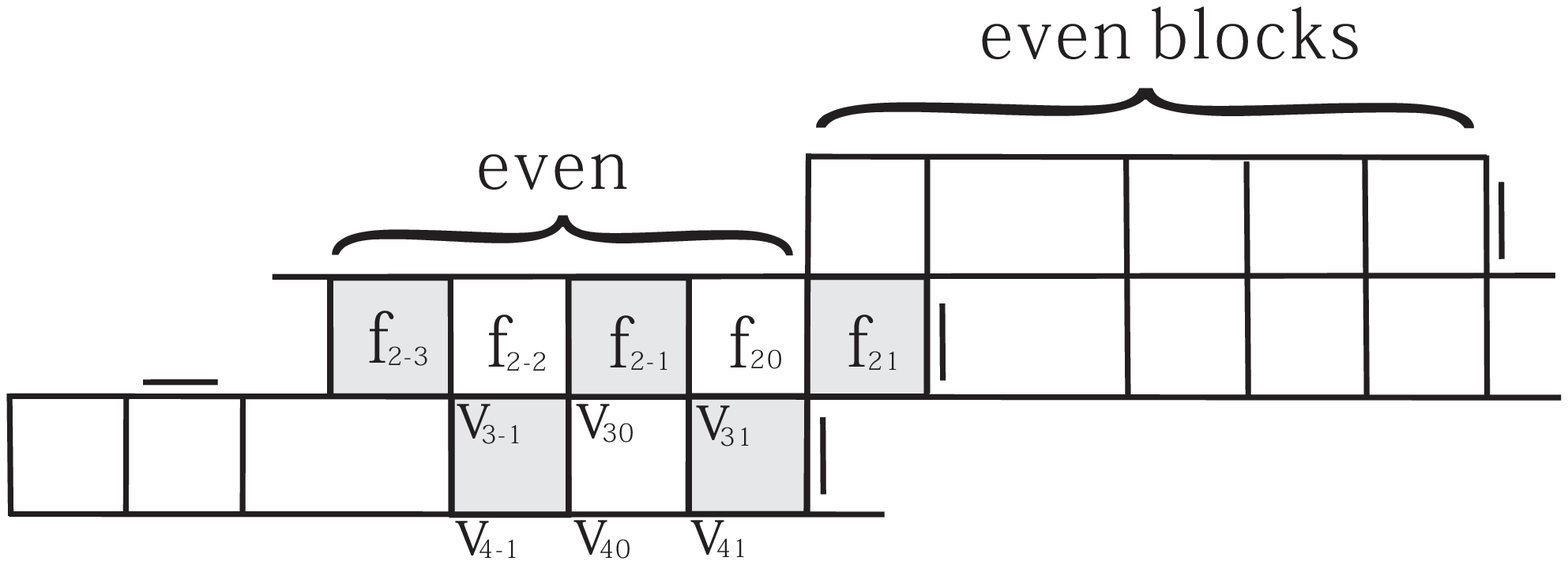}
\par {\small Fig.5.}
\end{figure}

\vskip 2mm
\begin{figure}[htbp]
\centering
\includegraphics[height=4cm]{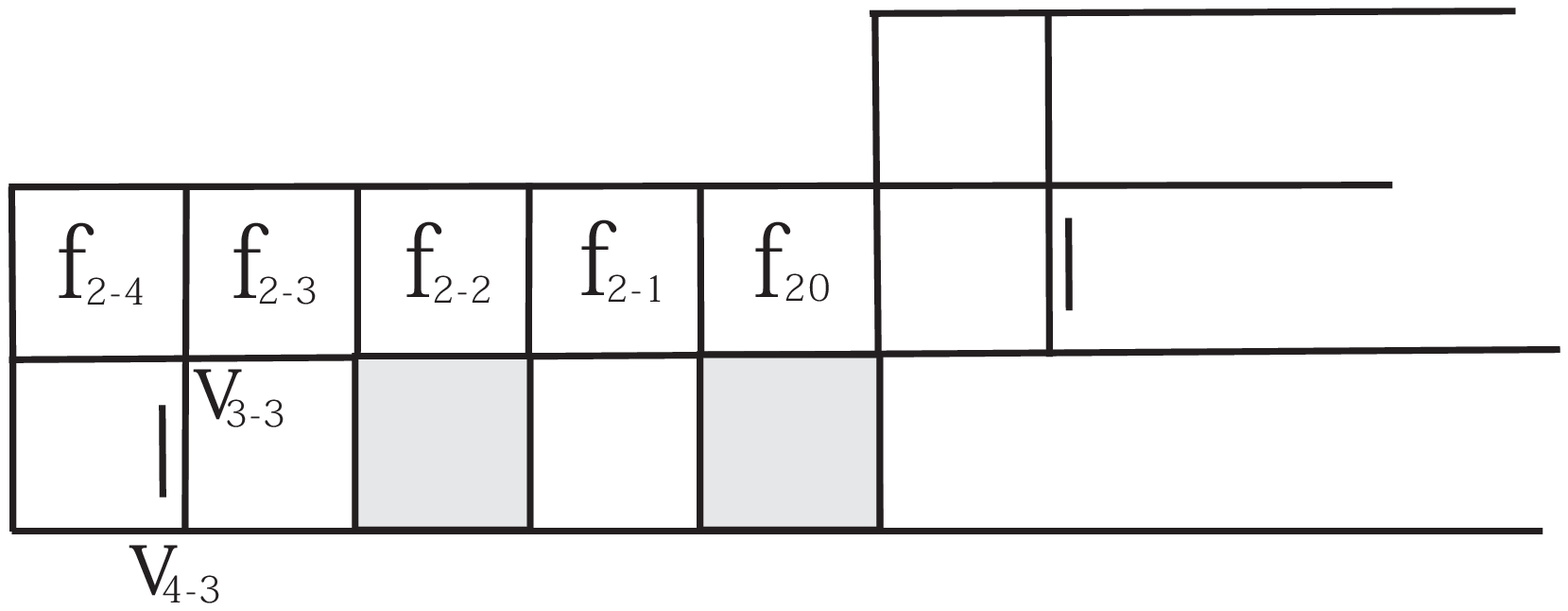}
\par {\small Fig.6.}
\end{figure}

\vskip 2mm
\begin{figure}[htbp]
\centering\centering
\includegraphics[height=6cm]{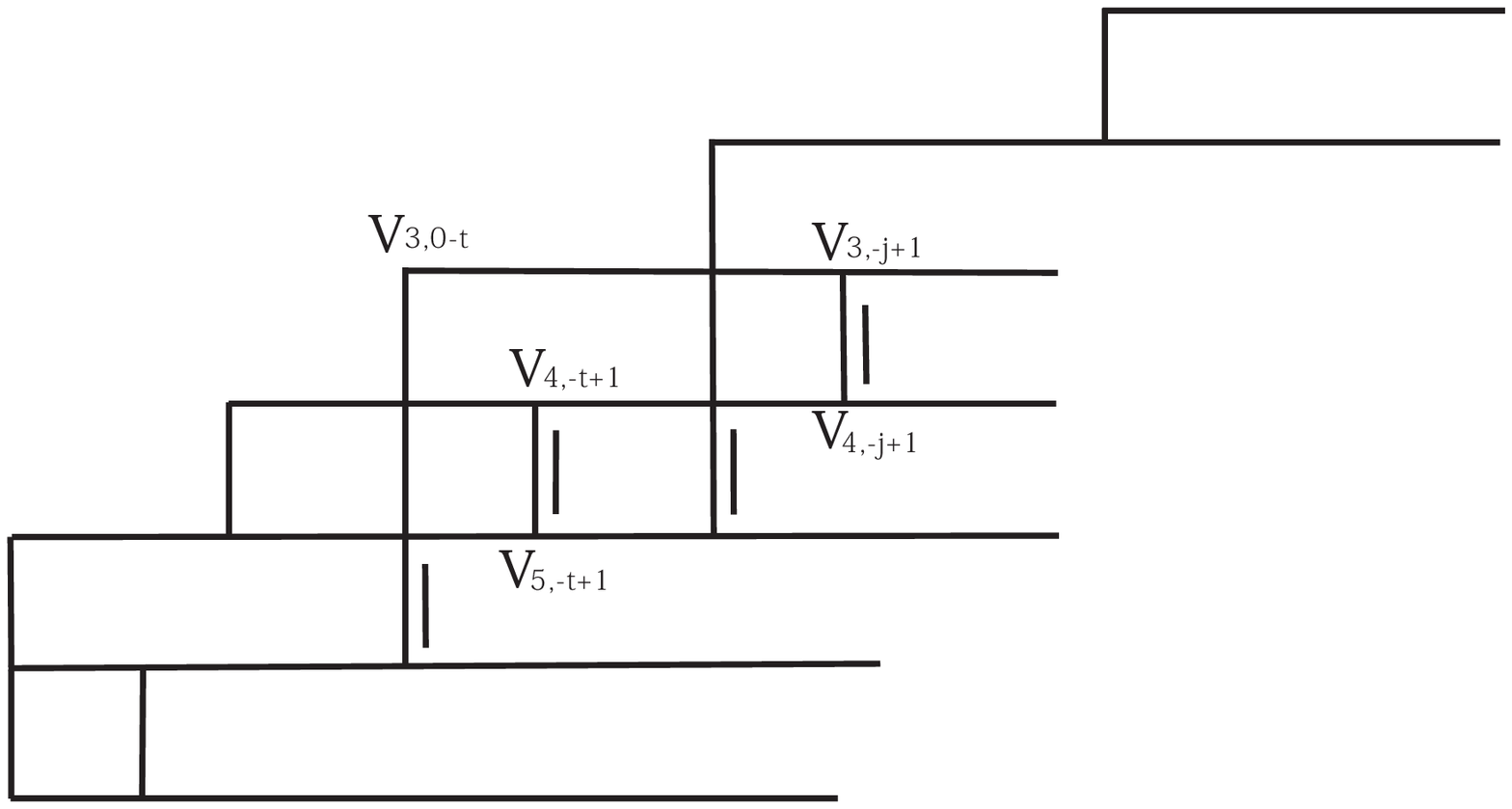}
\par {\small Fig.7.}
\end{figure}

In this case, we know that the edge $v_{4,-j}v_{5,-j}$ is either in the $M(G^*)$ or in $K'_{G^*}$, and alone the same line of reasoning as for the second row of $G^*$, we could show that either we could get a larger Clar cover of $G$ or $v_{4,-t+1}v_{5,-t+1}$ is in the perfect matching $M(G^*)$. The same argument terminates until the row containing the left most bottom block is encountered. i.e. the row containing the face $f_{x,-y}$ and there is no faces with $f_{x,-w}$ where $w \geq y$. See Fig.7 for detail. In this case, we find a larger Clar cover which has more tetragons than $C(G)$, a contradiction.
Consequently $G-C'$ has a unique
1-factor.$\Box$
\vskip 3mm

Now Let $M$ be a perfect matching of $G$ such that $f(G,M) = F(G)$. According to Theorem
2 in [20], there exist $F(G)$ disjoint $M$-alternating cycles, say as $C_1, C_2, . . . , C_{F(G)}$. Let $C^*$ denote
the subgraph of $G$ such that the outer boundary of $C*$ is the $M$-alternating cycle $C$.
From the proof of Theorems 9  in [20], it can be seen that, if the above theorem holds then $F(G)\leq C(G)$. Combines with Theorem 1,  the following conclusions hold.
\begin{theorem}  Let $G$ be an elementary polyomino. Then
$F(G)=C(G)$.
\end{theorem}

The following result follows immediately from Theorem 3 and Theorem 6.

\begin{theorem}  Let $G$ be an elementary polyomino. Then the maximum forcing number of $G$ can be computed in
polynomial time.
\end{theorem}

\section{The complexity of computing  the maximum forcing number of a non-elementary polyomino and hexagonal system}
\hspace*{8mm}

In this section we turn to the case of non- elementary polyomino and hexagonal system.

 In order to find the Clar number and maximum forcing number of a non-elementary polyomino and hexagonal system,
 we consider the decomposition of a non-elementary polyomino and hexagonal system
 into a number of elementary components.

 In {\cite{Zhang6}}, F.Zhang et al. developed an $O(n^2)$ algorithm to decompose a hexagonal system into a number of regions consisting of fixed bonds and a number of elementary components. Next, H.Zhang and F.Zhang{\cite{Zhang7}} also provided an algorithm to  decompose a polyomino into a number of regions consisting of fixed bonds and a number of elementary components. In \cite{Zhang8}, F.Zhang and H.Zhang generalized the above results. Let $G$ be a bipartite graph  whose vertices can be divided into two disjoint sets $U$ and $V$ (that is, $U$ and $V$ are each independent sets) such that every edge connects a vertex in $U$ to one in $V$. Colored each vertex in $U$ red and each vertex in $V$ blue. Let $M$ be a prefect matching of a bipartite graph $G$.  Orient all the edges belonging to $M$ toward the red vertices and orient the other edges of $G$ toward the blue vertices, the resulting digraph is called the orientation of $G$ and denoted by $D$. Let $e$ be an edge of $G$.   They showed that there is  an algorithm of $O(|E|+|V|)$ complexity to determine all elementary components and the fixed bonds of a bipartite graph $G$. The result is in the following:

\begin{theorem}[{\cite{Zhang8}}]  Let $G$ be a bipartite graph with a perfect matching $M$. A subgraph $H$ of $G$ is an elementary of $G$ iff the corresponding orientation of $H$ is a strongly connected component of $D$. And  there is an algorithm of $O(|E|+|V|)$ complexity to decompose $G$ into a number of regions consisting of fixed bonds and a number of elementary components.
\end{theorem}

Since any non-elementary polyomino (respectively, hexagonal system) with perfect matching can be composed into a number of elementary components and fixed bonds, then  the maximum forcing number of the original non-elementary polyomino (respectively, hexagonal system) with perfect matching is equal to the sum of the maximum forcing number of those elementary  components, and the Clar number of the original non-elementary polyomino(hexagonal systems) with perfect matchings is equal to the sum of the Clar number of those elementary  components. Thus we have the following results.

\begin{theorem}  Let $G$ be a polyomino or a hexagonal system with perfect matchings. Then
$F(G)=C(G)$.
\end{theorem}

Since the complexity of  decomposition
a non-elementary polyomino (respectively, hexagonal system) with perfect matchings into a number of elementary components and fixed bonds is $O(|E|+|V|)$, the number of elementary components of a non-elementary polyomino(respectively, hexagonal system) are no more than  $\frac{|V|}{2}$ and by Theorem 8, the maximum forcing number of every elementary components can be computed in polynomial time, then the maximum forcing number of a non-elementary polyomino or hexagonal system with perfect matching can be computed at $O(|E|+|V|)$ plus  polynomial time,  which are also polynomial time. Thus, we obtain the following results.

\begin{theorem}  Let $G$ be a polyomino(hexagonal systems) with perfect matchings. Then the maximum forcing number of $G$ can be computed in
polynomial time.
\end{theorem}

\end{document}